\documentclass[12pt]{article}
\usepackage{amssymb,amsmath,latexsym,color}
\usepackage{amscd,amsthm}
\newtheorem{example}[equation]{Example}
\newtheorem{theorem}[equation]{Theorem}

\newtheorem{lemma}[equation]{Lemma}
\newtheorem{proposition}[equation]{Proposition}
 
\newtheorem{remark}[equation]{Remark}

\numberwithin{equation}{section}

\newcommand{\fG}{\mathfrak{G}}

\newcommand{\ot}{\otimes}

\newcommand{\Z}{\mathbb{Z}}

\newcommand{\fh}{{\mathfrak h}}

\newcommand{\gs}{\sigma}

\newcommand{\End}{\hbox{End}}

\newcommand{\Hom}{{\rm Hom}}
\newcommand{\id}{{\rm id}}

\newcommand\Inn{\text{\rm Inn}}

\newcommand{\Tr}{\mathrm{Tr}}

\newcommand{\fsl}{\mathfrak{sl}}

\title{Jordan algebras and weight modules}

\author{Michael Lau$\hbox{}^{1}$\thanks{Corresponding author}\ \  and \ Olivier Mathieu$\hbox{}^{2}$\vspace{0.3cm}
\vspace{0.3cm}\\
$\hbox{\ \,}^1${\small Universit\'e Laval}\\ {\small D\'epartement de math\'ematiques et de statistique}\\{\small Qu\'ebec, QC,
Canada G1V 0A6}\\ {\small Michael.Lau@mat.ulaval.ca}\\
\\
$\hbox{}^{2 }${\small Institut Camille Jordan}\\ {\small UMR 5028 du CNRS, Universit\'e Claude Bernard Lyon 1}\\
{\small SUSTech International Center for Mathematics, Shenzhen, China}\\{\small 69622 Villeurbanne Cedex, France}\\ {\small mathieu@math.univ-lyon1.fr}\vspace{0.1cm}}

\date{}
\begin{document}
\maketitle
\setcounter{section}{0}

\begin{small}
  \noindent
  {\bf Abstract:} We consider bounded weight modules for the universal central extension $\fsl_2(J)$ of the Tits-Kantor-Koecher algebra of a unital Jordan algebra $J$.  Universal objects called Weyl modules are introduced and studied, and a combinatorial dominance criterion is given for analogues of highest weights.

  Specializing $J$ to the free Jordan algebra $J(r)$ of rank $r$, the category $\mathcal{C}^{fin}$ of finite-dimensional $\mathbb{Z}$-graded $\fsl_2(J)$-modules shares many properties with the representation theory of algebraic groups.  Using a deep result of Zelmanov, we show that this subcategory admits Weyl modules.  By analogy, we conjecture that $\mathcal{C}^{fin}$ is a highest weight category.  The resulting homological properties would then imply cohomological vanishing results previously conjectured as a way of determining graded dimensions of free Jordan algebras.

\bigskip


\bigskip

\end{small}
\maketitle

\section{Introduction}
Let $J$ be a unital Jordan algebra over an algebraically closed field $k$ of characteristic $0$.  Tits defined a Lie algebra
structure on the space $(\fsl_2(k)\otimes_kJ)\oplus \Inn\,J$, where
$\Inn\,J$ is the Lie algebra of inner derivations \cite{T}.
This construction was later generalized by  Kantor \cite{Kan}
and Koecher \cite{Koe}  and is now called the {\it Tits-Kantor-Koecher algebra} and  denoted by $TKK(J)$. 
The Lie algebra $TKK(J)$ is perfect and admits a universal central extension $\fsl_2(J)$ described \cite{AG}.  See also \cite{ABG,neher}.

With the exception of $r=1$ and $r=2$, the structure of the
free Jordan algebras $J(r)$  is unknown. However, it was proved in 
\cite{KM} that its structure is determined by the $\fsl_2(k)$-invariants of
$H_*(\fsl_2(J(r)))$. The following conjecture was provided.

\bigskip

\noindent
{\bf Conjecture A \cite[Conj.~3]{KM}.} {\em $H_n(\fsl_2(J(r)))^{\fsl_2(k)}=0,$ for all $n>0$.}

\bigskip
 
\noindent
Verification of this conjecture would give a recursive method to compute the dimensions of the graded components of $J(r)$.  In the present paper, we will 
interpret this  conjecture
in terms of representation theory.

Let $\{e, f,h\}$ be the standard basis of $\fsl_2(k)$, and let $h(a)=h\otimes a\in \fsl_2(J)$ for all $a\in J$.  For any $\fsl_2(J)$-module $M$ and integer $j$, let $M_j$ be the weight space $M_j=\{m\in M\,:\, h(1).m=j\,m\}.$  The module $M$ is said to be {\it bounded of level} $n$ if 
$$M=\bigoplus_{-n\leq j\leq n} M_j$$
for some nonnegative integer $n$ with $M_n\neq 0$.

A vector space $V$ endowed with a linear map
$\rho:J\to\End(V)$ is called a {\it $J$-space}
if it satisfies

\begin{enumerate}
\item[{\rm (J1)}] $[\rho(a),\rho(a^2)]=0$,
\item[{\rm (J2)}] $[[\rho(a),\rho(b)],\rho(c)]=  4\rho(\partial_{a,b}\,c)$,
\end{enumerate}

\noindent for all $a,b,c\in J$, where 
$\partial_{a,b}(c)=a(cb)-(ac)b$. When $\rho(1)$ acts by multiplication by $n$, $V$ is called
a {\it $J$-space of level $n$}. For any 
bounded $\fsl_2(J)$-module $M$ of level $n$, the weight space $M_n$ has a structure of
$J$-space of level $n$, where $\rho(a)$ is the action of $h(a)$. The $J$-space $V$ is said to be {\it dominant of level $n$}
if $V=M_n$ for some bounded $\fsl_2(J)$-module $M$ of level $n$. 

Our first main result characterizes 
dominant $J$-spaces of level $n$. For any
partition $\gs=(\gs_1,\gs_2,\ldots ,\gs_m)$ of $n+1$, we write $|C_\sigma|$ for the cardinality of the corresponding conjugacy class in $S_{n+1}$ and $\hbox{sgn}(\gs)$
for its  signature. We write $\rho_\sigma(a)$ for the expression $\rho(a^{\gs_1})\rho(a^{\gs_2})\cdots\rho(a^{\gs_m})$ for all $a\in J$. The following result appears as Theorem \ref{bounded-quotient-thm}.

\bigskip

\noindent
{\bf Theorem B.} {\em
Let $(V,\rho)$ be a $J$-space of level $n$.  Then 
$V$ is dominant if and only if
it satisfies the following   condition

$$\sum_{\gs\ \vdash\ n+1} \hbox{sgn}(\gs)\,|C_\gs|\,\rho_\gs(a)=0.$$}

\bigskip

For any dominant $J$-space of level $n$,
the {\it Weyl module} $\Delta(M)$ is the bounded
$\fsl_2(J)$-module of level $n$
defined by the following universal property: 
$$\Hom_{\fsl_2(J)}(\Delta(V),M)=\Hom_J(V,M_n),$$
for all bounded $\fsl_2(J)$-modules $M$ of level $n$, where homomorphisms of $J$-spaces are defined to be linear maps commuting with the action of $J$.

\bigskip

Assume now that $J=\oplus_{n\geq 0}\,J_n$ is a finitely generated $\mathbb{Z}_+$-graded unital Jordan algebra with $J_0=k1$.  The grading on $J$ clearly induces a grading on the Lie algebra $\fsl_2(J)$.  We show that the category $\mathcal{C}^{fin}$ of finite-dimensional $\mathbb{Z}$-graded $\fsl_2(J)$-modules admits a Weyl module for each dominant $J$-space.

\bigskip

\noindent
    {\bf Theorem C.} {\em For any $\Z$-graded finite-dimensional dominant $J$-space $V$ of level $n$, the Weyl module $\Delta(V)$ is finite-dimensional.}

\bigskip

\noindent
Theorem C is nontrivial and uses a deep result of Zelmanov on nil Jordan algebras.  It appears as Theorem \ref{fin-dim-thm} in the paper, and shows that    the category
$\mathcal{C}^{fin}$ shares many properties with categories of representations of
reductive algebraic groups in positive characteristic.  This leads to the following conjecture.

\bigskip

\noindent
{\bf Conjecture D.} {\em The category $\mathcal{C}^{fin}$ is a highest weight category, in the sense of Cline, Parshall, and Scott.}

\bigskip

\noindent
A proof of Conjecture D would also settle Conjecture A of \cite{KM}.

\bigskip

\noindent
    {\bf Theorem E.}  {\em Let $r\geq 1$. If $\mathcal {C}^{fin}$ is a highest weight category, then Conjecture A holds.}

\bigskip
    
    \noindent
    Theorem E appears below as Theorem \ref{thm-3}.

\section{Bounded weight modules}

Let $J$ be a unital Jordan algebra over an algebraically closed field $k$ of characteristic zero.  All vector spaces, algebras, and tensor products will be taken over $k$.  For every $a,b\in J$, let $L_a:\ J\rightarrow J$ be the multiplication operator $L_a(b)=ab$.  Write $\partial_{a,b}=[L_a,L_b]$, and let $\kappa$ be one half the Killing form on $\mathfrak{sl}_2(k)$.  We write $\Inn\,J$ for the set $\{\partial_{a,b}\,:\,a,b\in J\}$ of {\em inner derivations} of $J$.  The element $x\ot a$ in the vector space $\fsl_2(k)\ot J$ will be denoted by $x(a)$.  We fix a standard basis $\{h,e,f\}$ of $\fsl_2(k)$ with $[h,e]=2e$, $[h,f]=-2f$, and $[e,f]=h$.

\subsection{Tits construction}

In his 1962 paper, Tits defined a Lie algebra structure on the space 
$$TKK(J):=
\fsl_2(k)\otimes J\oplus \Inn\,J,$$
with Lie bracket
\begin{enumerate} 
\item{\rm (T1)} 
$[x(a),y(b)]=[x,y](ab)+ \kappa(x,y)\partial_{a,b}$

 \item{\rm (T2)} $[\partial, x(a)]=x(\partial\,a)$,
\end{enumerate}
where $x(a)=x\ot a$ 
\noindent for any $x,y\in\fsl_2$, $\partial\in\Inn\,J$, and $a,b\in J$.  This construction was later generalized to Jordan pairs and triple systems by Kantor and Koecher, and $TKK(J)$ is known as the {\em Tits-Kantor-Koecher (TKK) algebra}.

    \subsection{Tits-Allison-Gao construction}

The TKK algebra is perfect, that is,  
$TKK(J)=[TKK(J),TKK(J)]$, so $TKK(J)$ admits a universal central extension, which we denote by $\fsl_2(J)$.  This Lie algebra was nicely described in the 1996 paper of Allison and Gao \cite{AG} in the context of universal coverings of the Steinberg unitary Lie algebras $\mathfrak{stu}_n(J)$ for $n\geq 3$.  The case where $n=3$ corresponds to $TKK(J)$.  See also \cite{ABG} for equivalent formulas written in terms of $\fsl_2(k)$.

As a vector space,
$$\mathfrak{sl}_2(J)=\left(\mathfrak{sl}_2(K)\otimes J\right)\oplus\{J,J\},$$

\noindent where $\{J,J\}=(\bigwedge^2 J)/\mathcal{S}$ and
$\mathcal{S}=\hbox{Span}\{a\wedge a^2\mid \ a
\in J\}$.
For any $a,b\in J$, we write 
$\{a,b\}$ for the image of  $a\wedge b$ in $\{J,J\}$.
The bracket on $\mathfrak{sl}_2(J)$ is given by

\begin{enumerate}
\item[{\rm (R1)}]$[x(a),y(b)]=[x,y](ab)+
\kappa(x,y)\{a,b\}$

\item[{\rm (R2)}]$[\{a,b\},x(c)]=x(\partial_{a,b}\,c)$

\item[{\rm (R3)}]$[\{a,b\},\{c,d\}]=\{\partial_{a,b}\,c,d\}+\{c,\partial_{a,b}\,d \}$.
\end{enumerate}

\noindent for all $a,b,c,d\in J$ and $x,y\in\mathfrak{sl}_2(k)$.
 It is a bit tricky to show that 
(R3) is skew-symmetric \cite{AG}.

There is an obvious Lie algebra epimorphism
$\fsl_2(J)\to TKK(J)$ which is the identity on
$\fsl_2\otimes J$ and sends the symbol $\{a,b\}$ to $\partial_{a,b}$.  When the Jordan algebra $J$ is associative, we have $\{J,J\}=HC_1(J)$, and the construction specializes to results of Kassel and Loday\cite{KL}.\footnote{The notation in \cite{KL} differs slightly from the modern conventions--Kassel-Loday write $HC_2$ for what is now denoted as $HC_1$.}

\subsection{The short grading of  $\fsl_2(J)$}

The Lie algebra $\mathfrak{G}:=\mathfrak{sl}_2(J)$ decomposes with respect to the adjoint action of $\frac12 h(1)$ as 
$\mathfrak{G}=\mathfrak{G}_{-1}\oplus\mathfrak{G}_{0}\oplus\mathfrak{G}_{1}$, where

\hskip3.7cm $\mathfrak{G}_{-1}=f\otimes J$,

\hskip3.7cm $\mathfrak{G}_{0}=h\otimes J\oplus \{J,J\}$,
and

\hskip3.7cm $\mathfrak{G}_{1}=e\otimes J$.

\noindent This decomposition is a root grading in the sense of Berman-Moody \cite{BM}, and is called the {\it short grading} of 
$\mathcal{G}$.  In fact, every root-graded Lie algebra of type $A_1$ has a universal central extension isomorphic to $\fsl_2(J)$ for some unital Jordan algebra $J$.  See \cite{neher} or \cite{ABG} for details.

\subsection{Bounded modules}

For any $\fsl_2(J)$-module $V$ and any $k\in\Z$,
let

\centerline{$V_k=\{v\in V\,:\,  h(1)v=k\,v\}$.}

\noindent An $\fsl_2(J)$-module $V$ is said to be a 
          {\it bounded weight module of level $\ell$} if $V=\oplus_{-\ell\leq m\leq \ell}\, V_m,$ with $V_\ell\neq 0$.

\subsection{$J$-spaces}

\noindent Recall that a vector space $M$ endowed with a linear map
$\rho:\ J\rightarrow\End(M)$ is called
 a {\it $J$-space} if $\rho$ satisfies

\begin{enumerate}
\item[{\rm (J1)}] $[\rho(a),\rho(a^2)]=0$,
\item[{\rm (J2)}] $[[\rho(a),\rho(b)],\rho(c)]=  4\rho(\partial_{a,b}\,c)$.
\end{enumerate}

\noindent
A $J$-space is said to be of {\em level n} if $\rho(h(1))=n\,\id$.

%Any $\fG_0$-module structure on a vector space $V$ can 
%been seen as a Lie algebra homomorphism
%$\tilde{\rho}:\fG_0\to \End(V)$. In what follows, we will identify $J$ with the subspace $h\otimes J\subset\fG_0$. 

\begin{lemma}\label{J-space-lemma} 
Let $(M,\tilde{\rho})$ be a representation of the Lie algebra $\fG_0$. Then the map
$\rho:a\mapsto\tilde{\rho}(h(a))$ determines a $J$-space structure on $M$.

Conversely, if $M$ is a $J$-space, then there is
a unique
$\fG_0$-module structure $(M,\tilde{\rho})$ such that 
$\tilde{\rho}(h(a))=\rho(a)$ for any $a\in J$.
\end{lemma}

\noindent
{\bf Proof.} Let $(M,\tilde{\rho})$ be a $\fG_0$-module.  For $a\in J$, set $\rho(a)=\tilde{\rho}(a)$.  Since $[h(a),h(a^2)]=4\{a,a^2\}=0$, it follows 
that $[\rho(a),\rho(a^2)]=0$, proving
{\rm (J1)}. Let $a,b,c\in J$.
We have $[[h(a),h(b)],h(c)]=4h(\partial_{a,b}\,c)$,
and therefore 
 $[[\rho(a),\rho(b)],\rho(c)]=  4\rho(\partial_{a,b}\,c)$, proving {\rm (J2)}.
 
 Conversely, assume that $M$ is a $J$-space.
 It is clear that $\fG_0$ is  generated by
 the vector space $h\otimes J$ and defined by the relations
 
 \begin{enumerate}
 
\item[{\rm(H1)}] $[h(a),h(a^2)]=0$, and
 
\item[{\rm(H2)}]  $[[h(a),h(b)],h(c)]=4h(\partial_{a,b}\,c)$,
 
\end{enumerate}
 
\noindent for any $a,b,c\in J$.
Therefore there is a unique structure
$\tilde{\rho}$ of
$\fG_0$-module on $M$ such that
$\tilde{\rho}(h(a))=\rho(a)$\qed

\bigskip

A linear map $\gs:\ J\rightarrow \End_k(M)$ is a {\em Jordan birepresentation} (and $M$ is a {\em Jordan bimodule}) if the semidirect product $(a,m)(b,n)=(ab,bm+an)$ gives a Jordan algebra structure to the vector space direct sum $J\oplus M$.  This condition is equivalent to the conditions

 \begin{enumerate}
 
\item[{\rm(M1)}] $[\gs(a),\gs(a^2)]=0$, and
 
\item[{\rm(M2)}]  $\gs(a^2b)+2\gs(a)\gs(b)\gs(a)=2\gs(ab)\gs(a)+\gs(a^2)\gs(b)$,
 
\end{enumerate}
 for all $a,b\in J$.  It follows from \cite[II.9(47')]{Jacobson} that (M1) and (M2) imply that $\rho=2\gs$ satisfies (J1) and (J2).  The converse is not true, however.  If $(M,\rho)$ is a $J$-space of level $n$, then setting $a=b=1\in J$ in (M2), we see that the only possible values for $n=\rho(1)$ are $0,1$ or $2$ (Peirce decomposition), if we wish to induce a $J$-bimodule structure on $M$ with action $\gs=\frac12\rho$.  But as we will see in Example 2.10, there exist $J$-spaces of any level.

\subsection{Dominant $J$-spaces}

Any $J$-space $M$ of level $n$ can be induced to a generalized Verma module
$$V(M)=\mathcal{U}(\fG)\ot_{\mathcal{U}(\fG_0\oplus\fG_1)}M,$$
where $\fG_1$ acts as zero on $M$.  The main result of this section will be an analysis of which $J$-spaces $M$ of level $n$ determine $\fG$-modules $V(M)$ with bounded quotients $V(M)/X$, such that the composition of natural maps
$$M\hookrightarrow V(M)\rightarrow V(M)/X$$
is an injection of $\fG_0\oplus\fG_1$-modules.  Such a quotient $V(M)/X$ will be called a {\em bounded $M$-quotient} of $V(M)$ of level $n$, and in this case, $M$ is said to be {\em dominant}.  Note that if $V(M)/X$ is a bounded $M$-quotient of level $n$, then $n$ is a nonnegative integer and the weight $n$ subspace $(V(M)/X)_n$ is precisely $M$.

\begin{proposition}\label{bounded-prop}
  Let $M$ be a $J$-space of nonnegative integer level $n$, with $\fG_0$-action given by $\rho:\ J\rightarrow\hbox{End}_k(M)$.  Then the generalized Verma module $V(M)$ has a bounded $M$-quotient if and only if
$$e(1)^{n+1}f(a)^{n+1}m=0,$$
for all $a\in J$ and $m\in M$.
\end{proposition}

\noindent
{\bf Proof.}  If $V(M)$ has a nonzero bounded $M$-quotient $V(M)/X$, then $M=(V(M)/X)_n=V(M)_n$, so $X_n=0$.  By $\mathfrak{sl}_2$-theory, we see that $(V(M)/X)_m=0$ for all $m<-n$, so $X_m=V(M)_m$ for all $m<-n$.  In particular, $f(a)^{n+1}m\in V(M)_{-n-2}=X_{-n-2}\subseteq X$ for all $m\in M$, so $e(1)^{n+1}f(a)^{n+1}m\in X$ for all $a\in J$ and $m\in M$.  But then $e(1)^{n+1}f(a)^{n+1}m\in V(M)_n\cap X=0$.

Conversely, suppose that $e(1)^{n+1}f(a)^{n+1}m=0$ for all $a\in J$ and $m\in M$.  Let $Z\subset V(M)$ be a $\fG$-submodule which is maximal with respect to the property that $Z_n=0$.  Linearisation of the relation $e(1)^{n+1}f(a)^{n+1}m=0$ gives
$$e(1)^{n+1}f(b_1)\cdots f(b_{n+1})m=0,$$
for all $b_1,\ldots,b_{n+1}\in J$ and $m\in M$.  Then for any $a_1,\ldots,a_{n+1}\in J$, we have
$$h(a_1)\ldots h(a_{n+1})e(1)^{n+1}f(b_1)\cdots f(b_{n+1})m=0,$$ from which it follows that
$$e(a_1)\cdots e(a_{n+1})f(b_1)\cdots f(b_{n+1})m=0,$$
for all $a_i,b_j\in J$ and $m\in M$.  In particular, this shows that $f(b_1)\cdots f(b_{n+1})m\in Z$, and $V(M)_k\subseteq Z$ for all $k\leq -n-2$.  Therefore, $V(M)_k=Z_k$ for all $k<-n$, and $V(M)/Z$ is a bounded $M$-quotient.\qed

\bigskip

We now introduce some notation.  Let $\gs=(\gs_1,\ldots,\gs_m)$ be a partition of $n+1$, that is,
 $$\gs_1\geq\cdots\geq\gs_m\geq 1,\ \hbox{for some}\ m\geq 1,\ \hbox{where}\ \gs_1+\cdots+\gs_m=n+1.$$
We write $|C_\sigma|$ for the cardinality of the conjugacy class $C_\sigma$ of permutations in the symmetric group $S_{n+1}$ with cycle structure $\gs$.  The sign of these permutations will be denoted by $\hbox{sgn}(\gs)$, and we write $\rho_\sigma(a)$ for the expression $\rho(a^{\gs_1})\rho(a^{\gs_2})\cdots\rho(a^{\gs_m})$ for all $a\in J$ and $\gs\ \vdash\ n+1$.  It follows easily from Condition (C1) that $[\rho(a^i),\rho(a^j)]=0$ for all $i,j$, so this product is independent of the order of the factors.

The Newton polynomials $N_\ell(x)=x_1^\ell+\cdots +x_n^\ell$ in $n$ indeterminates $x_1,\ldots,x_n$ generate the ring of symmetric polynomials $k[x_1,\ldots,x_n]^{S_n}$.  We write
$$N_\gs(x)=N_{\gs_1}(x)N_{\gs_2}(x)\cdots N_{\gs_m}(x)$$
for each partition $\gs=(\gs_1,\ldots,\gs_m)$ of $n+1$.  The space $P_{n+1}\subset k[x_1,\ldots,x_n]^{S_n}$ of symmetric polynomials of total degree $n+1$ is clearly of dimension $p(n+1)-1$, where $p(n+1)$ is the number of partitions of $n+1$.  The Newton polynomials $N_1(x),\ldots, N_n(x)$ are algebraically independent, so
$$\{N_\gs(x)\ :\ \gs\ \vdash\ n+1\ \hbox{such that}\ \gs\neq(1,1,\ldots,1)\}$$
is a basis for $P_{n+1}$, and the set
$$\{N_\gs(x)\ :\ \gs\ \vdash\ n+1\}$$
has exactly one linear dependence relation, up to scalar multiple.  We include an amusing representation-theoretic argument below, that we have not seen elsewhere in the literature.

\begin{proposition}\label{unique-lin-combo-prop}
Up to scalar multiple, the unique linear dependence relation on the set $\{N_\gs(x)\ :\ \gs\ \vdash\ n+1\}$ is $\displaystyle{\sum_{\gs\ \vdash\ n+1} \hbox{sgn}(\gs)|C_\gs|N_\gs(x)=0}.$
\end{proposition}

\noindent
{\bf Proof.}  By the Frobenius character formula,
$$N_\gs(x)=\sum_{\lambda}\chi_\lambda(\gs)S_\lambda,$$
where the sum is taken over all partitions $\lambda\neq (1,1,\ldots,1)$ of $n+1$.  Here $\chi_\lambda(\gs)$ is the character (evaluated at any permutation of cycle structure $\gs$) of the Specht module associated with $\lambda$, and $S_\lambda$ is the Schur polynomial associated to $\lambda$.  See \cite{FH} for details.

Since $\hbox{sgn}$ is the character of the sign representation, the Specht module associated to $(1,\ldots,1)$, we see that
\begin{align*}
\sum_{\gs\ \vdash\ n+1}\hbox{sgn}(\gs)|C_\gs|N_\gs(x)&=\sum_{\gs\in S_{n+1}}\hbox{sgn}(\gs)N_\gs(x)\\
&=\sum_{\gs \in S_{n+1}}\chi_{(1,\ldots,1)}(\gs)\sum_{\lambda\neq(1,\ldots,1)}\chi_\lambda(\gs)S_\lambda\\
&=\sum_{\lambda\neq(1,\ldots,1)}S_\lambda\sum_{\gs\in S_{n+1}}\chi_{(1,\ldots,1)}(\gs)\chi_\lambda(\gs),
\end{align*}
and the inner product
$$(\chi_{(1,\ldots,1)},\chi_\lambda)=\sum_{\gs\in S_{n+1}}\chi_{(1,\ldots,1)}(\gs)\chi_\lambda(\gs)$$
is $0$ whenever $\lambda\neq(1,\ldots,1)$, by the orthogonality relations.\qed

\bigskip

\begin{remark}{\em 
  Proposition \ref{unique-lin-combo-prop} can also be proved with a more standard combinatorial argument:  Let $V$ be a vector space of dimension $n$, with basis
$\{e_1,\ldots,e_n\}$. Let $h:V\to V$ be defined by
$h(e_i)=x_i e_i$, and let  
$h^{\otimes (n+1)}:V^{\otimes (n+1)}
\to V^{\otimes (n+1)}$ be the induced endomorphism.
For $\sigma\in S_{n+1}$, we have
$$\Tr (h^{\otimes (n+1)}\circ\sigma)=N_\sigma(x).$$
Since $\sum_{\sigma\in S_{n+1}}\hbox{sgn}(\gs)\sigma$ acts as zero on $V^{\otimes (n+1)}$, the dependence relation follows.}
  \end{remark}

\bigskip

The following well-known formula, originally due to Garland \cite{Garland78} and reinterpreted by Chari and Pressley \cite{CP01}, will be used to prove boundedness conditions.

\begin{lemma} \label{gen-fn-formula} Let $p:\ \mathcal{U}(\fG)\rightarrow \mathcal{U}(\fG_{-1})\mathcal{U}(\fG_0)$ be the projection relative to the (vector space) decomposition
$\mathcal{U}(\fG)=\mathcal{U}(\fG)\fG_1\oplus\mathcal{U}(\fG_{-1})\mathcal{U}(\fG_0).$  Then $p(e(1)^rf(a)^{n+1})$ is the coefficient of $u^{n+1}$ in the generating function
$$\frac{(-1)^r r!(n+1)!}{(n+1-r)!}\left(\sum_{s=1}^\infty f(a^s)u^s\right)^{n+1-r}\exp\left(-\sum_{t=1}^\infty \frac{h(a^t)}{t}u^t\right).$$
\end{lemma}\qed

\begin{theorem}\label{bounded-quotient-thm}
Let $(M,\rho)$ be a $J$-space of level $n$.  Then the following conditions are equivalent:
\begin{enumerate}
\item[{\rm (1)}] $M$ is dominant.
\item[{\rm (2)}] $e(1)^{n+1}f(a)^{n+1}m=0$ for all $a\in J$ and $m\in M$.
\item[{\rm (3)}] $\displaystyle{\sum_{\gs\ \vdash\ n+1} \hbox{sgn}(\gs)|C_\gs|\rho_\gs(a)=0}.$
\end{enumerate}
\end{theorem}

\noindent
{\bf Proof.}  By Proposition \ref{bounded-prop}, Conditions (1) and (2) are equivalent, so we need only prove that (2) and (3) are equivalent.

%Now assume that Condition (2) holds.
Since $e(1)^{n+1}f(a)^{n+1}$ is homogeneous of degree $0$ with respect to the grading induced by $\hbox{ad}\,(h\ot 1)$, we see that its action on any highest weight vector $m$ is given by the action of its projection $p$ on the subspace $\mathcal{U}(\fG_0)$ with respect to the decomposition $\mathcal{U}(\fG)_0=\mathcal{U}(\fG_0)\oplus\left(\mathcal{U}(\fG)\fG_1\cap\mathcal{U}(\fG)_0\right)$ of the space $\mathcal{U}(\fG)_0$ of degree $0$ elements of $\mathcal{U}(\fG)$.

By Lemma \ref{gen-fn-formula}, $p(e(1)^{n+1}f(a)^{n+1})$ is the coefficient of $u^{n+1}$ in the generating series
$$(-1)^{n+1}(n+1)!(n+1)!\hbox{exp}\left(-\sum_{k=1}^\infty\frac{h(a^k)}{k}u^k\right).$$
Computing directly, the coefficient of $u^{n+1}$ in $\hbox{exp}\left(-\sum_{k=1}^\infty\frac{h(a^k)}{k}u^k\right)$ is
$$\sum_{\gs\vdash n+1}(-1)^{r_\gs}\frac{h_{\gs}(a)}{\left(\prod_{i=1}^{r_\gs}\gs_i\right)\left(\prod_{j=1}^{m_\gs}a_j!\right)},$$
where $\gs_i$ is the length of the $i$th row of the Young frame $T_\gs$ associated to the partition $\gs$, $r_\gs$ is the number of rows of $T_\gs$, $m_\gs$ is the number of columns of $T_\gs$, $a_j$ is the number of rows of length $j$ in $T_\gs$, and $h_\sigma(a)=h(a^{\sigma_1})h(a^{\sigma_2})\cdots h(a^{\sigma_{r_\sigma}})$.  If $\hbox{odd}(r_\gs)$ (respectively, $\hbox{even}(r_\gs)$ is the number of odd-length (respectively, even-length) rows of $T_\gs$, we see that $(-1)^{n+1}=(-1)^{\hbox{odd}(r_\gs)}$, so
$$(-1)^{n+1}(-1)^{r_\gs}=(-1)^{n+1}(-1)^{\hbox{odd}(r_\gs)}(-1)^{\hbox{even}(r_\gs)}=(-1)^{\hbox{even}(r_\gs)}=\hbox{sgn}(\gs).$$
By elementary counting arguments,
$$|C_\gs|=\frac{(n+1)!}{\left(\prod_{i=1}^{r_\gs}\gs_i\right)\left(\prod_{j=1}^{m_\gs}a_j!\right)}.$$
See \cite[Proposition 1.1.1]{Sagan01}, for instance.  The projection of $e(1)^{n+1}f(a)^{n+1}$ on $\mathcal{U}(\fG_0)$ is thus  $\displaystyle{(n+1)!\sum_{\gs\ \vdash\ n+1} \hbox{sgn}(\gs)|C_\gs|h_\gs(a)},$ so Conditions (2) and (3) are equivalent.\qed

% which acts as zero on every $J$-space of level $n$.

%Applying this observation to the Jordan subalgebra\footnote{In fact, this algebra is actually associative.}  $A\subseteq J$ generated by $1$ and $a$, we see that there is a nontrivial linear combination
%$$\sum_{\gs\ \vdash\ n+1}b_\gs h_\gs(a)$$
%which acts as zero on every $A$-space of level $n$.  In particular, $k[x_1,\ldots,x_n]$ is an $A$-space of level $n$ if we let $\rho(a^\ell)$ act as multiplication by the Newton polynomial $N_\ell(x)$ for each $\ell\geq 0$.  By Proposition \ref{unique-lin-combo-prop}, we now see that
%$$b_\gs=\hbox{sgn}(\gs)|C_\gs|,$$
%for all $\gs\ \vdash\ n+1$.  Any $J$-space $M$ of level $n$ is also an $A$-space of level $n$, so by uniqueness of the relation in Proposition \ref{unique-lin-combo-prop}, the same expression $\sum\hbox{sgn}(\gs)|C_\gs|h_\gs(a)$ acts as zero on $M$.  That is,
%$$\sum_{\gs\ \vdash\ n+1}\hbox{sgn}(\gs)|C_\gs|\rho_\gs(a)=0.$$

%That (3) implies (2) is clear, since we have now shown that $e(a)^{n+1}f(1)^{n+1}m=\sum_{\gs\ \vdash\ n+1}\hbox{sgn}(\gs)|C_\gs|\rho_\gs(a)m,$ by the uniqueness of the relation in Proposition \ref{unique-lin-combo-prop}.\qed

%Since $e(a)m=0$ for all $a\in J$ and $m\in M$, it follows from the Poincar\'e-Birkhoff-Witt theorem that $e(a)^{n+1}f(1)^{n+1}m$ can be expressed as some linear combination of $\{\rho_\gs(a)m\ :\ \gs\ \vdash\ n+1\}$.

\bigskip

\begin{example}\label{bounded-J-space-example}{\em By Theorem \ref{bounded-quotient-thm}, any dominant $J$-space $M$ of level $0$ is trivial, in the sense that $\rho:\ J\rightarrow \hbox{End}_k(M)$ is the zero map and $M$, equipped with the trivial $\fG$-action, is the unique bounded $M$-quotient of $V(M)$.}
\end{example}

\begin{example}{\em 
    Dominant $J$-spaces $(M,\rho)$ of level $1$ satisfy $\rho(a^2)=\rho(a)^2$ for all $a\in J$, so
    \begin{equation}\label{assoc-spec}
    \rho(ab)=\frac12(\rho(a)\rho(b)+\rho(b)\rho(a)),
\end{equation}
    for all $a,b\in J$ by linearization.  %Condition (\ref{assoc-spec}) implies Conditions (C1) and (C2) of Lemma \ref{J-space-lemma}, so 
    Dominant $J$-spaces $(M,\rho)$ of level $1$ are thus precisely {\em associative specializations}, Jordan algebra homomorphisms $\rho$ from $J$ to special Jordan algebras of linear operators on a vector space $M$.}
\end{example}

\begin{example}{\em 
    For levels higher than $2$, dominant $J$-spaces are never Jordan bimodules.  See the discussion after Lemma \ref{J-space-lemma} for details.  Many such $J$-spaces exist.  For example, it follows immediately from Proposition \ref{unique-lin-combo-prop} and Theorem \ref{bounded-quotient-thm}, that the map
    \begin{align}
      \rho:\ k[t]&\rightarrow \hbox{End}_k\left(k[x_1,\ldots,x_n]^{S_n}\right)\\
      t^\ell&\mapsto N_\ell(x)\nonumber
\end{align}
    defines a dominant $J$-space of level $n$ for the (associative) Jordan algebra $J=k[t]$, an example we will consider in more detail in Section 3.

}\end{example}

%\bigskip
%\noindent
%{\it Remark} There exists a notion of 
%a Jordan module investigated in \cite{Jacobson}.
%Let $V$ be a Jordan module, and let
%$\sigma:J\to\End(V)$ be the structural map.
%By definition, we have

%\centerline{$[\sigma(a),\sigma(a^2)]=0$ and}

%\centerline{$[[\sigma(a),\sigma(b)]\sigma(c)]=
%\sigma(\partial_{a,b}\,c)$,}

%\noindent for any $a,b,c\in J$. It follows 
%that $\rho:=2\sigma$ defines a structure of $J$-space
%on $V$. Conversely, if $\rho$ is a structure of $J$-space on $V$, $1/2 \rho$ is usually not a Jordan module. We will see later on the precise connection between these two notions.

\section{Weyl modules and highest weight categories}

%For any $J$-space $(M,\rho)$, the {\em multiplication algebra} $\mathcal{M}(M)$ is the associative (but not necessarily commutative) subalgebra generated by $\{\rho(a)\ :\ a\in J\}\subseteq\hbox{End}_k(M)$.  A $J$-space $(M,\rho)$ is {\em free of rank $r$} if there exist $m_1,\ldots, m_r\in M$ such that $M=\sum_{i=1}^r\mathcal{M}(M)m_i$, and there is no nontrivial linear combination $\sum_{i=1}^rp_im_i=0$ with $p_i\in\mathcal{M}(M)$ for all $i$.

Let $n$ be a nonnegative integer.  The categories $\mathcal{C}^b(M)$ of bounded weight modules attached to bounded $J$-spaces $M$ of level $n$ admit universal objects $\Delta(M)$, called {\em Weyl modules}.  Every bounded $M$-quotient of level $n$ is a homomorphic image of $\Delta(M)$, and it is clear that
$$\Delta(M)=V(M)\ \left/\ \right.\mathcal{U}(\fG)\sum_{\ell<-n}V(M)_\ell\ \ =\ \ \bigoplus_{\ell=0}^n\Delta(M)_{n-2\ell},$$
where $\Delta(M)_{n-2\ell}$ is the vector subspace of weight $n-2\ell$.  Identifying
$$\fG_{-1}=\{f\ot a\ :\ a\in J\}\subset\fG=(\mathfrak{sl}_2(k)\ot J)\oplus\{J,J\}$$ with $J$, the weight space $V(M)_{n-2\ell}$ identifies with the vector space $S^\ell J\ot M$ for $\ell=0,\ldots,n$.  %In this section, we consider Weyl modules $\Delta_n(M)$, concentrating on the case where $M$ is a free $J$-space for a free Jordan algebra $J$ of small rank.

\subsection{Weyl modules for finite dimensional dominant $J$-spaces}
Let $J=\bigoplus_{\ell=0}^\infty J_\ell$ be a finitely generated $\mathbb{Z}_+$-graded unital Jordan algebra with $J_0=k1$.  Let $M=\bigoplus_{\ell=0}^\infty M_\ell$ be a $\mathbb{Z}$-graded dominant $J$-space of level $n$.  We now prove one of our main results, that the category $\mathcal{C}^{fin}$ of finite-dimensional $\mathbb{Z}$-graded $\fsl_2(J)$-modules contains its Weyl modules.

\begin{theorem}\label{fin-dim-thm}  Let $M$ be a $\mathbb{Z}$-graded dominant $J$-space of level $n$, for a $Z_+$-graded and finitely generated Jordan algebra $J$ with $J_0=k1$.  Then the Weyl module $\Delta(M)$ is finite dimensional if and only if $M$ is finite dimensional.
\end{theorem}

\noindent
{\bf Proof.}  If $\Delta(M)$ is finite dimensional, then $M\subseteq \Delta(M)$ is clearly also finite dimensional.  Conversely, assume that $M$ is finite dimensional.  Up to a possible shift in grading, we may assume that $M$ is $\mathbb{Z}_+$-graded.   Let $N$ be the largest nonnegative integer for which the graded component $M_N$ is nonzero.  Let $a\in J$ be a homogeneous element with $\deg a>N$, and let $v\in M$.  By Lemma \ref{gen-fn-formula}, $e(1)^nf(a)^{n+1}v$ is the coefficient of $u^{n+1}$ in the formal series 
$$(-1)^nn!(n+1)!\sum_{s=1}^\infty f(a^s)u^s\exp\left(-\sum_{t=1}^\infty \frac{h(a^t)}{t}u^t\right)v.$$
By degree considerations, $h(a^t)v=0$ for all $t\geq 1$, so
$$e(1)^nf(a)^{n+1}v=(-1)^nn!(n+1)!f(a^{n+1})v,$$
and $f(a^{n+1})v=0$ as an element of $\Delta(M)$.

In particular, $f\left(b^{(N+1)(n+1)}\right)M=0$ for all $b$ in the (non-unital) Jordan subalgebra $J^+=\bigoplus_{\ell=1}^\infty J_\ell\subset J$.  Let $I=\{x\in J^+\,:\,f(x)M=0\}$.  For all $x\in J^+$, $y\in I$, and $m\in M$,
\begin{align*}
0&=h(x)f(y)m\\
&=f(y)h(x)m-2f(xy)m\\
&=-2f(xy)m
\end{align*}
since $h(x)M\subseteq M$ and $f(y)M=0$.  Therefore, $xy\in I$ and $I$ is an ideal of $J^+$.

Since $b^{(N+1)(n+1)}\in I$ for all $b\in J^+$, the Jordan algebra $J^+/I$ is nil of bounded index, hence locally nilpotent by a result of Zelmanov \cite{Zelmanov79}.  But $J$, and thus $J^+/I$, is finitely generated, so $J^+/I$ is nilpotent and there exists $N'>0$ such that every product (in any association) of $N'$ elements of $J^+$ is in $I$.  The (finitely many) generators of $J$ may be chosen to be homogeneous and of positive degree at most $r$ for some $r>0$.  In particular, $J_s\subseteq I$ for all $s\geq rN'$.  That is, $f(a)M=0$ for all $a\in J$ with $\deg a\geq rN'$.

The weight space $\Delta(M)_{n-2\ell}$ is spanned by monomials of the form 
$$f(a_1)\cdots f(a_\ell)w$$ with $a_1,\ldots,a_\ell\in J$ and $w\in M.$  
Since the $f(a_i)$ commute with each other, the set 
$$\{f(a_1)\cdots f(a_\ell)w\,:\,w\in M \hbox{\ and\ } a_1,\ldots,a_\ell\in J \hbox{\ with\ }\deg a_i<rN'\hbox{\ for all\ }i\}$$ 
already spans $\Delta(M)_{n-2\ell}$.  As $M$ and $J_i$ are finite dimensional for all $i$, it now follows that $\dim \Delta(M)_{n-2\ell}<\infty$ for all $\ell$, and the Weyl module $\Delta(M)=\bigoplus_{\ell=0}^n\Delta(M)_{n-2\ell}$ is also finite dimensional.\qed

\subsection{Highest weight categories and character formulas for free Jordan algebras}

Cline, Parshall, and Scott \cite{cps88} introduced the notion of {\em highest weight category} as a unifying theme in representation theory, modelled after highest weight representations of semisimple algebraic groups and their Lie algebras.  Their definition requires labelling simple objects by a poset $\Lambda$, and the existence of enough injectives, as well as costandard objects labelled by the same index set as the simples and satisfying various axioms.  Given the similarities between the category $\mathcal{C}^{fin}$ of finite-dimensional $\mathbb{Z}$-graded $\fsl_2(J(r))$-modules and the representation theory of reductive algebraic groups in positive characteristic, we conjecture that $\mathcal{C}^{fin}$ is a highest weight category, with the Weyl modules and their duals (twisted by the Cartan involution) as the standard and costandard objects, respectively.  In a highest weight category, the higher ext-groups $\hbox{Ext}^i(\Delta(\lambda),\nabla(\mu))=0$ for all $i>0$ and $\lambda,\mu\in\Lambda$.  If $\mathcal{C}^{fin}$ is indeed a highest weight category as conjectured above, the vanishing of higher ext-groups would, in fact, settle the main conjecture of \cite{KM} and thus describe the graded dimensions of the free Jordan algebras $J(r)$.

%A positive answer to this conjecture would, in fact, imply the vanishing of the relative cohomology $H^i(\fsl_2(J),\fsl_2(k))$ for $i>0$, settling the main conjecture in \cite{KM}, and thus determining character formulas for the graded components of the free Jordan algebra $J(r)$.

\bigskip

\begin{theorem} \label{thm-3}
  If $\mathcal{C}^{fin}$ is a highest weight category with Weyl modules and their duals as its standard and costandard objects, then $H_i(\fsl_2(J(r)))$ contains no nonzero trivial $\fsl_2(k)$-modules for $i>0$.
\end{theorem}

\noindent
    {\bf Proof.} Let $J$ be the free unital Jordan algebra $J(r)$ on $r$ generators, and suppose that $\mathcal{C}^{fin}$ is a highest weight category as in the hypotheses of the theorem.  As noted above, in a highest weight category, $\hbox{Ext}^i(\Delta(\lambda),\nabla(\mu))$ vanishes for all $i>0$ and indices $\lambda,\mu$ of simples, where $\Delta(\lambda)$ and $\nabla(\mu)$ are the corresponding standard and costandard objects.  In $\mathcal{C}^{fin}$, the Weyl and dual Weyl modules corresponding to the trivial $1$-dimensional $\fsl_2(J)$-module $k$ are themselves $1$-dimensional, so
    $$\hbox{Ext}_{\mathcal{C}^{fin}}^i(k,k)=0\hbox{\quad for all\ }i>0.$$
    But $\hbox{Ext}_{\mathcal{C}^{fin}}^i(k,k)=H^i(\fsl_2(J))$, and the cohomology ring
    $$H^*(\fsl_2(J))= H^*(\fsl_2(k))\ot H^*(\fsl_2(J),\fsl_2(k)).$$
    As $H^0(\fsl_2(J))=H^0(\fsl_2(k))=k$, we see that $H^*(sl_2(J))=k$ and the relative cohomology $H^i(\fsl_2(J),\fsl_2(k))=0$, for all $i>0$.  The result now follows from the universal coefficient theorem and the interpretation of the relative cohomology as the $\fsl_2(k)$-invariants in $H^i(\fsl_2(J))$.\qed

\subsection{Example: Weyl modules for free Jordan algebras of rank 1}

For any Jordan algebra $J$ with unit $1$ and $n\in \mathbb{Z}_+$, let
$$T(J)=k1\oplus J\oplus \left(J^{\ot 2}\right)\oplus \left(J^{\ot 3}\right)\oplus\cdots$$
be its tensor algebra, and let $I\subseteq T(J)$ be the two-sided ideal generated by the relations

\begin{align}
&1-n 1,\\
  &a\ot a^2-a^2\ot a,\\
&a\ot b\ot c+c\ot b\ot a-b\ot a\ot c-c\ot a\ot b+b(ac)-a(bc),\\
&\sum_{\sigma\vdash n+1}\hbox{sgn}(\gs)|C_\gs|T_\gs(a),\label{3.4}
\end{align}
for all $a,b,c\in J$, where $T_\gs(a)=a^{\gs_1}\ot a^{\gs_2}\ot\cdots\ot a^{\gs_{m}}$ for all partitions $\gs=(\gs_1,\ldots ,\gs_{m})\vdash n+1$.  The associative algebra $\mathcal{U}_n(J)=T(J)/I$ is called the {\em universal $J$-space envelope of level $n$}.  There is a unique associative algebra homomorphism $\check{\rho}:\ T(J)\rightarrow\hbox{End}_k(M)$ extending the action $\rho:\ J\rightarrow\hbox{End}_k(M)$ of any $J$-space $(M,\rho)$ of level $n$, and in light of Lemma \ref{J-space-lemma} and Theorem \ref{bounded-quotient-thm}, the map $\check{\rho}$ descends to the quotient $\mathcal{U}_n(J)$.  By construction, dominant $J$-spaces of level $n$ and left $\mathcal{U}_n(J)$-modules are equivalent notions, and a $J$-space $(M,\rho)$ of level $n$ is said to be {\em free of rank $r$} if $(M,\check{\rho})$ is a free $\mathcal{U}_n(J)$-module.

Let $F=\mathcal{U}_n(J)$ be the universal $J$-space envelope of level $n$ for a unital Jordan algebra $J$.  If $J$ is finitely generated as a Jordan algebra, then $F$ is finitely generated as an associative algebra, by Relation (\ref{3.4}).  For example, if $J=k[t]$ is the free Jordan algebra of rank $1$, then $F=k[x_1,\ldots,x_n]^{S_n}$ is the algebra of symmetric polynomials, where $t^\ell$ corresponds to the Newton polynomial $N_\ell(x)=x_1^\ell+\cdots+x_n^\ell\in F$.  If $J$ is free of rank $m$, then $\mathcal{U}_0(J)=k$ and $\mathcal{U}_1(J)$ is the quotient of the free associative algebra in $m$ generators by the ideal generated by the relation $a\ot b\ot c+c\ot b\ot a=b\ot a\ot c+c\ot a\ot b$ for all $a,b,c\in J$.
%In Example \ref{bounded-J-space-example}, we saw that
%\begin{align*}
%      \rho:\ k[t]&\rightarrow \hbox{End}_k\left(F\right)\\
%      t^m&\mapsto N_m(x)
%\end{align*}
%defines a $J$-space structure on the vector space $F=k[x_1,\ldots,x_n]^{S_n}$ of symmetric polynomials, where $J=k[t]$ is the free Jordan algebra of rank $1$.  The Newton polynomials $N_0(x),\ldots ,N_n(x)$ generated the ring $k[x_1,\ldots,x_n]^{S_n}$ of symmetric polynomials as a $k$-algebra, so the image of $\check{\rho}:\ T(J)\rightarrow\hbox{End}_k(F)$ is precisely $k[x_1,\ldots,x_n]^{S_n}$.  The algebraic relations among Newton polynomials in $n$ variables are generated by the dependence relation of Proposition \ref{unique-lin-combo-prop}, so the kernel of $\check{\rho}$ is in the ideal $I$.  That is, $F$ is a free $\mathcal{U}_n(J)$-module of rank $1$, generated by the element $1\in F$.

Let $L$ be the two-dimensional simple $\mathfrak{sl}_2(k)$-module.  The Jordan algebra $J=k[t]$ is commutative and associative, and it is easy to see that $\{J,J\}=0$ and the TKK algebra $\fG=\mathfrak{sl}_2(J)=\mathfrak{sl}_2(k)\ot J$ is centrally closed.  The space $L[t]=L\ot k[t]$ is obviously a $\fG$-module, where
$$(x\ot p(t)).(v\ot q(t))=xv\ot p(t)q(t),$$
for all $x\in\mathfrak{sl}_2(k)$, $v\in L$, and $p(t),q(t)\in J$.  This gives a $\fG$-module structure on the space $S^n(L[t])\subset T(L[t])$ of homogeneous symmetric tensors of degree $n$.

\begin{proposition}Let $F=k[x_1,\ldots,x_n]^{S_n}$ be the rank $1$ free $\mathcal{U}_n(k[t])$-module.  Then the Weyl module $\Delta_n(F)$ is isomorphic to $S^n(L[t])$.  \end{proposition}
{\bf Proof.}  Let $v\in L$ be a nonzero vector of weight $1$ with respect to the action of $h\in \mathfrak{sl}_2(k)$.  There is a natural injection $\iota:\ F\longrightarrow S^n(L[t]),$ with
$$  \iota:\ \sum_{\gs\in S_n}x_{\gs(1)}^{a_{\gs(1)}}\cdots x_{\gs(n)}^{a_{\gs(n)}}\longmapsto\sum_{\gs\in S_n}(v\ot t^{a_{\gs(1)}})\ot\cdots\ot(v\ot t^{a_{\gs(n)}}). $$
This maps extends uniquely to a $\fG$-module epimorphism
\begin{align*}
  V(F)&\rightarrow S^n(L[t])\\
  u.p&\mapsto u.\iota(p),\hbox{\quad\quad for all\ }u\in U(\fG)\hbox{\ and\ }p\in F,
  \end{align*}
with kernel $\sum_{\ell<-n}V(F)_\ell$.\qed

\begin{remark}{\em
    In fact, for every prime Jordan algebra $J$ and every $n\geq 2$, there is a dominant $J$-space of level $n$, on which the Lie algebra $\fG_0(J)=(\fh\ot J)\oplus\{J,J\}$ acts faithfully.  If $J$ is special, then there is a faithful associative specialization $\rho:\ J\rightarrow \hbox{End}_k(M)$, and $M$ is a $J$-space of level $1$.  The faithfulness of the extension $\tilde{\rho}:\ \fG_0\rightarrow \hbox{End}_k(M)$ on $\{J,J\}$ follows immediately from the assumption that $J$ is prime.  We can then take the $n$-fold tensor product of $M$ to obtain a faithful $J$-space of level $n$.

    If $J$ is the Albert algebra $\mathbb{A}$, then we can construct a faithful $\fG_0$-module of level $n$ as a tensor product of copies of the level $2$ and level $3$ representations of the Albert algebra, obtained from representations of the exceptional Lie algebra $E_6$, viewed as the subalgebra $(k\,h\ot\mathbb{A})\oplus\{\mathbb{A},\mathbb{A}\}$ of the Lie algebra $\fsl_2(\mathbb{A})$.

    This observation is clearly not true for arbitrary (non-prime) Jordan algebras.  For example, the Lie algebra $\fG_0(k[t,t^{-1}])$ has a nontrivial centre that acts as $0$ on all bounded modules.}
\end{remark}

%\begin{proposition}
%Combinatorial description of weight spaces of $\Delta(M)$ when $M$ is $1$-dimensional trivial module and $J=k[t]$.
%\end{proposition}
%{\bf Proof.} 

%\begin{theorem}
%Description of weight spaces of $\Delta(M)$ for $J=k[t]$ and $M$ arbitrary.
%\end{theorem}
%{\bf Proof.} 

%\subsection{Example: Weyl modules for free Jordan algebras of rank $2$}
%Connection to generation of s-identities

%\section{Extensions in the category $C^b$}
%Relative cohomology and Ext in the bounded and unbounded category.
%Connection to growth conjecture. 

%\bigskip

%\noindent
%Michael Lau, D\'epartement de math\'ematiques et de statistique, Universit\'e Laval, Qu\'ebec, Canada G1V0A6, Michael.Lau@mat.ulaval.ca

%\bigskip

%\noindent
%Olivier Mathieu, Institut Camille Jordan, UMR 5028 du CNRS, Universit\'e Claude Bernard Lyon 1, 69622 Villeurbanne Cedex, France, mathieu@math.univ-lyon1.fr

\end{document}